\documentclass{amsart}
\usepackage{amsmath,amsthm,amssymb,IMjournal}
\usepackage{graphicx}

\begin{document}

\newcommand{\bff}{\boldsymbol{f}}
\newcommand{\bfg}{\boldsymbol{g}}

\newcommand{\Erelper}{\mathrm{RelErr}(\mathrm{period})}
\newcommand{\ErelL}{\mathrm{RelErr}(L^2)}

\newcommand{\fmax}{f^{\mathrm{max}}}
\newcommand{\fmin}{f^{\mathrm{min}}}
\newcommand{\rd}{\mathrm{d}}
\newcommand{\porig}{p_\mathrm{orig}}
\newcommand{\papprox}{p_\mathrm{approx}}

\newcommand{\ttg}[1]{\textrm{(\ref{#1})}}

\newtheorem{theorem}{Theorem}[section]

\numberwithin{equation}{section}

\title{Accurate reduction of a model of circadian rhythms
by delayed quasi steady state assumptions}

\author{\|Tom\'a\v{s} |Vejchodsk\'y|, Oxford}

\rec {September 9, 2013}


\abstract 
Quasi steady state assumptions are often used to simplify complex systems
of ordinary differential equations in modelling of biochemical processes.
The simplified system is designed to have the same qualitative properties
as the original system and to have a small number of variables. This enables to
use the stability and bifurcation analysis to reveal a deeper structure in
the dynamics of the original system. This contribution shows that introducing
delays to quasi steady state assumptions yields a simplified system that
accurately agrees with the original system not only qualitatively but also
quantitatively. We derive the proper size of the delays for a particular model
of circadian rhythms and present numerical results showing the accuracy
of this approach.
\endabstract



\thanks
   The research leading to these results has received funding from the
   People Programme (Marie Curie Actions) of the European Union's Seventh
   Framework Programme (FP7/2007-2013) under REA grant agreement no.~328008.
   Further, I acknowledge the support of RVO 67985840.
   Finally, let me thank R.\ Erban and P.K.\ Maini for introducing me this
   topic and for fruitful discussions.
\endthanks

\section{Introduction}\label{sec1}


Model reduction is a crucial technique in large biochemical systems,
because it enables to employ analytical and numerical methods
to reveal detailed structure of the kinetics \cite{CotVejErb:2013,ErbChaKevVej:2009}.
As an example, we consider 
a theoretical biochemical model of circadian rhythms described in \cite{Vilar:2002}.
Using the law of mass action \cite{Savageau:1969}, the kinetics
of this chemical system can be
described by a system of nine nonlinear ordinary differential equations (ODEs).
The authors of \cite{Vilar:2002} use various quasi steady state
assumptions to reduce the system to just two ODEs in such a way that the
reduced system has the same qualitative behaviour as the original system,
i.e. a periodic solution.
Then, they use the reduced system to perform the bifurcation and stochastic
analysis.

However, using parameters introduced in \cite{Vilar:2002}, the period
of the original system is about 25 hours while the period of the reduced
system is roughly 18 hours. Thus, the relative error in the period is
approximately 30\,\%. The error in the amplitude is even close to 100\,\%,
as we show in Table~\ref{ta:1} below.

In this contribution we study the quasi steady state assumptions in detail.
We use numerical quadrature to derive explicit formulas for
delays for approximated chemical species and reduce the original system of
nine chemical reactions to two delay ODEs. Some of the delays
depend on the state variables in a complicated way, which might be problematic
for the subsequent analysis, therefore we show that this dependence can be simplified.
Finally, numerical solutions show that periods
of the original and reduced system agree within 2\,\% relative error
and that the error in the amplitude decreases to about 20\,\%.


The following section introduces the model of circadian rhythms
and its quasi steady state reduction.
Section~\ref{se:justification} justifies the quasi steady state assumptions
and Section~\ref{se:derivation} derives the delayed quasi steady state assumption.
The accuracy of these approximations
is assessed in Section~\ref{se:numer} and final conclusions are drawn
in Section~\ref{se:conclusion}.

\section{Model of circadian rhythms}
\label{se:model}

\begin{figure}
\includegraphics[width=0.5\textwidth]{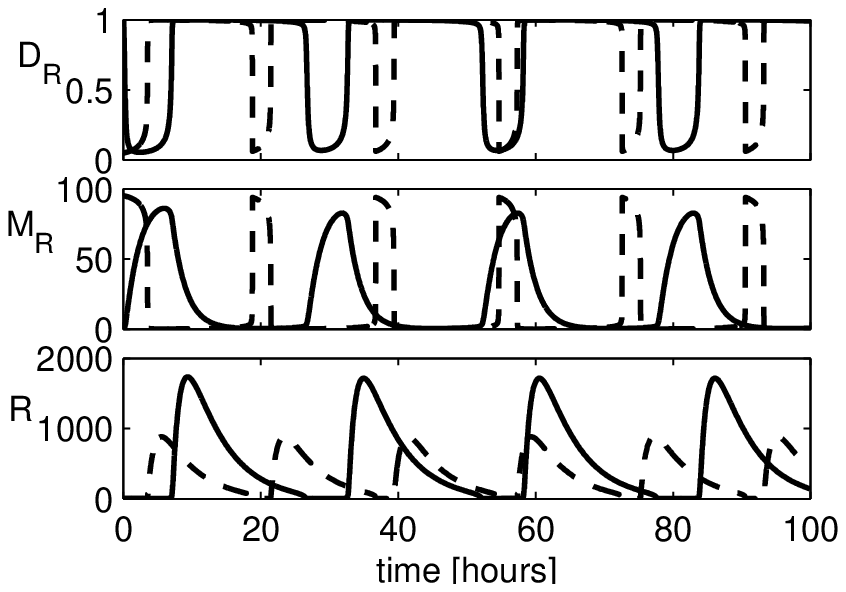}%
\includegraphics[width=0.5\textwidth]{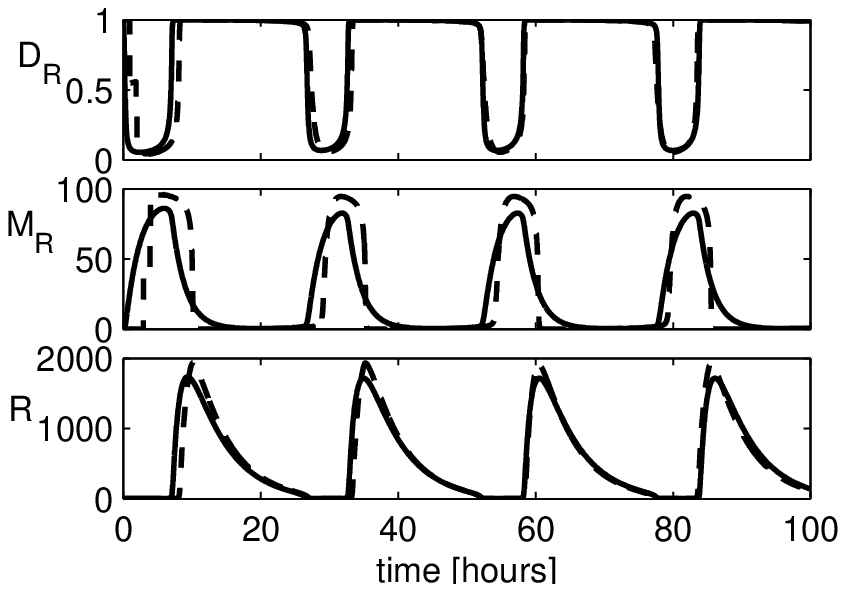}
\caption{\label{fi:1}
Components $D_R$, $M_R$, and $R$ of the solution 
of the original system \ttg{eq:DA}--\ttg{eq:C} (solid lines),
reduced system \ttg{eq:rR}--\ttg{eq:rC} (dashed lines, left panel),
and delay system \ttg{eq:dR}--\ttg{eq:dC} (dashed lines, right panel).
}
\end{figure}

Circadian rhythms are modelled in \cite{Vilar:2002} by the following
system of nine ODEs:
\begin{align}
\label{eq:DA}
  \rd  D_A / \rd t  &=
  \theta_A D_A' - \gamma_A  D_A  A,
\\ \label{eq:DAprime}
 \rd  D_A' / \rd t  &=
 -\theta_A  D_A' + \gamma_A  D_A  A,
\\ \label{eq:DR}
  \rd  D_R / \rd t  &=
  \theta_R D_R' - \gamma_R  D_R  A,
\\ \label{eq:DRprime}
  \rd  D_R' / \rd t  &=
  -\theta_R  D_R' + \gamma_R  D_R  A,
\\ \label{eq:MA}
  \rd  M_A / \rd t  &=
  \alpha_A'  D_A' + \alpha_A  D_A - \delta_{M_A}  M_A,
\\ \label{eq:MR}
  \rd  M_R / \rd t  &=
  \alpha_R' D_R' + \alpha_R  D_R - \delta_{M_R}  M_R,
\\ \label{eq:A}
  \rd  A / \rd t  &=
  \beta_A  M_A + \theta_A D_A' + \theta_R { D_R'}
    -  A(\gamma_A  D_A + \gamma_R  D_R + \gamma_C  R + \delta_A),
\\ \label{eq:R}
  \rd  R / \rd t  &=
  \beta_R  M_R - \gamma_C  A  R + \delta_A  C - \delta_R  R,
\\ \label{eq:C}
  \rd  C / \rd t  &=
  \gamma_C  A  R - \delta_A  C.
\end{align}
Here, the time variable is denoted by $t$, the capital letters stand for
the copy numbers of respective molecules that evolve in time and Greek letters
stand for the rate constants. Namely, $A=A(t)$, $R=R(t)$, $M_A=M_A(t)$,
$M_R=M_R(t)$, $D_A = D_A(t)$, $D_R = D_R(t)$ denote the numbers of molecules
of the activator, represor, their mRNA and genes, respectively. Functions
$D_A' = D_A'(t)$ and $D_R' = D_R'(t)$ stand for the number of molecules
of the activated forms of genes and $C=C(t)$ for the complex of $A$ and $R$.
Values for the rate constants are taken from \cite{Vilar:2002} as
\begin{multline}
\label{eq:param}
\alpha_A = 50\,\mathrm{h}^{-1},\ 
\alpha_A' = 500\,\mathrm{h}^{-1},\ 
\alpha_R = 0.01\,\mathrm{h}^{-1},\
\alpha_R' = 50\,\mathrm{h}^{-1},\
\beta_A = 50\,\mathrm{h}^{-1},\\
\beta_R = 5\,\mathrm{h}^{-1},\
\gamma_A = 1\,\mathrm{Mol}^{-1} \mathrm{h}^{-1},\
\gamma_R = 1\,\mathrm{Mol}^{-1} \mathrm{h}^{-1},\
\gamma_C = 2\,\mathrm{Mol}^{-1} \mathrm{h}^{-1},\
\delta_A = 1\,\mathrm{h}^{-1},\\
\delta_R = 0.2\,\mathrm{h}^{-1},\
\delta_{M_A} = 10\,\mathrm{h}^{-1},\
\delta_{M_R} = 0.5\,\mathrm{h}^{-1},\
\theta_A = 50\,\mathrm{h}^{-1},\
\theta_R = 100\,\mathrm{h}^{-1}.
\end{multline}
Notice that by Mol and h we understand the number of molecules and the hour.
The initial condition for system \ttg{eq:DA}--\ttg{eq:C} is considered as
\begin{equation}
\label{eq:icond}
D_A =  D_R = 1\,\mathrm{Mol},\quad
D_A' =  D_R' =  M_A =  M_R =  A =  R =  C = 0\,\mathrm{Mol}.
\end{equation}
Figure~\ref{fi:1} shows three components of the solution of system \ttg{eq:DA}--\ttg{eq:C} with
parameter values \ttg{eq:param} and initial condition \ttg{eq:icond}
as solid lines.

To reduce the system, let us first notice that
$\rd (D_A + D_A') / \rd t = 0$ and $\rd (D_R + D_R') / \rd t = 0$.
Thus, due to the initial condition we infer conservation relations
\begin{equation}
\label{eq:subst}
  D_A' = 1 - D_A \quad\text{and}\quad D_R' = 1 - D_R
\end{equation}
that enable to eliminate $D_A'$ and $D_R'$ from the system by simple substitution.
To simplify the system further we use so-called quasi steady state assumptions
\cite{Murray:2002,SegelSlemrod:1989}.

In general, the idea of quasi steady state assumptions is based on splitting the system
into slow and fast variables. The steady state of fast variables depends on
values of slow variables. If slow variables change, the steady states
change as well and the fast variables go quickly towards their new steady states.
Thus it is a reasonable approximation to consider the fast variables to be
effectively in their steady states. Of course, the quality of this approximation
depends on actual speeds of the dynamics of the slow and fast variables.

In case of system \ttg{eq:DA}--\ttg{eq:C}, we simply assume that $D_A$, $D_R$,
$M_A$, $M_R$, and $A$ are fast and stay in their steady states that may however
change with the values of the slow variables $R$ and $C$.
From \ttg{eq:DA}, \ttg{eq:DR}, \ttg{eq:MA}, and \ttg{eq:MR} with \ttg{eq:subst},
we easily obtain steady states for $D_A$, $D_R$, $M_A$ and $M_R$ as functions of $A$:
\begin{alignat}{2}
\label{eq:steadystates}
  D_A^s(A) &= \frac{\theta_A}{\theta_A + \gamma_A A},
  &\qquad 
  M_A^s(A) 
    &= \frac{\alpha_A'}{\delta_{M_A}} + \frac{\theta_A(\alpha_A - \alpha_A')}{\delta_{M_A}(\theta_A+\gamma_A A)},
\\ \nonumber
  D_R^s(A) &= \frac{\theta_R}{\theta_R + \gamma_R A},
  &\qquad
  M_R^s(A) 
    &= \frac{\alpha_R'}{\delta_{M_R}} + \frac{\theta_R(\alpha_R - \alpha_R')}{\delta_{M_R}(\theta_R+\gamma_R A)}.
\end{alignat}
Approximating $D_A$, $D_R$, $M_A$, and $M_R$ by their respective steady states
in \ttg{eq:A}, we can express the steady state of $A$ as the following function of $R$:
\begin{equation}
\label{eq:tA}
  \widetilde A^s(R) = \frac12(\alpha_A' \rho(R) - K_d)
  + \frac12\sqrt{(\alpha_A' \rho(R)-K_d)^2 + 4\alpha_A\rho(R)K_d},
\end{equation}
where $\rho(R) = \beta_A/ (\delta_{M_A} (\gamma_C R + \delta_A))$
and $K_d = \theta_A/\gamma_A$, see \cite{Vilar:2002}.
Using the approximation $A = \widetilde A^s(R)$, we may express steady
states \ttg{eq:steadystates} as functions of $R$ and reduce the original system
\ttg{eq:DA}--\ttg{eq:C} to just two ODEs for $R$ and $C$:
\begin{align}
\label{eq:rR}
  \rd  R / \rd t  &=
  \beta_R  M_R^s(\widetilde A^s(R)) - \gamma_C  \widetilde A^s(R)  R + \delta_A  C - \delta_R  R,
\\
\label{eq:rC}
  \rd  C / \rd t &= \gamma_C \widetilde A^s(R) R - \delta_A  C.
\end{align}
Figure~\ref{fi:1} (left panel) shows $R(t)$ as the solution of \ttg{eq:rR}--\ttg{eq:rC}
together with approximations of $D_R = D_R^s(\widetilde A^s(R))$ and
$M_R = M_R^s(\widetilde A^s(R))$
as dashed lines.

\section{Justification of quasi steady state assumptions}
\label{se:justification}

Let us justify the above described quasi steady state assumptions
on an illustrative example of equation \ttg{eq:MA}.
Using \ttg{eq:subst}, we can express equation \ttg{eq:MA} as
\begin{equation}
\label{eq:MA2}
  \rd  M_A(t) / \rd t  = \Phi(t)  - \delta_{M_A}  M_A(t),
  \quad \text{where }
  \Phi(t) = \alpha_A' + (\alpha_A - \alpha_A')  D_A(t).
\end{equation}
If the function $\Phi(t)$ was explicitly known, we could easily find an
expression for the solution $M_A(t)$ to \ttg{eq:MA2} with the initial condition
\ttg{eq:icond} as
\begin{equation}
  \label{eq:MAsol}
  M_A(t) = \int_0^t \Phi(\tau) \exp[\delta_{M_A}(\tau-t)] \,\rd\tau.
\end{equation}
To obtain a quasi steady state approximation of $M_A(t)$ we employ one-point
numerical quadrature in \ttg{eq:MAsol} and approximate
$M_A(t) \approx w_1 \Phi(t_1)$, where $t_1 = t$ is the quadrature point
and the corresponding quadrature weight $w_1$ is determined
such that the resulting rule integrates constant functions exactly:
$w_1 = \int_0^t \exp[\delta_{M_A}(\tau-t)]\,\rd\tau = (1 - \exp(-\delta_{M_A}t))/\delta_{M_A}$.
Since the exponential function $\exp(-\delta_{M_A}t)$ decays quickly to zero,
we may neglect it with respect to 1 for sufficiently large $t$. As a result
we approximate $w_1\approx 1/\delta_{M_A}$ and $M_A(t) \approx \Phi(t)/\delta_{M_A}$
which is exactly the quasi steady state approximation $M_A^s(A)$
from \ttg{eq:steadystates} provided $D_A = D_A^s(A)$.

\section{Derivation of delayed quasi steady state assumptions}
\label{se:derivation}

The reasoning from the previous section can be made more accurate, because
one-point quadrature rules have the capability to be exact even for all
linear functions.
More precisely, we consider a quadrature point
$t_2 = t - \tau_{M_A}$, a corresponding weight $w_2$,
and approximate the integral in \ttg{eq:MAsol} by
$w_2 \Phi(t_2)$. We find the particular values of $\tau_{M_A}$ and $w_2$
such that this quadrature rule is exact for all linear functions.

Any linear function can be expressed as
$\ell(\tau) = \alpha_1 \ell_1(\tau) + \alpha_2 \ell_2(\tau)$, where
$\ell_1(\tau) = (\tau + \tau_{M_A} - t)/\tau_{M_A}$ and
$\ell_2(\tau) = (t - \tau)/\tau_{M_A}$.
Clearly, $\ell$ is determined by its values $\alpha_1$ and $\alpha_2$
at points $t_1 = t$ and $t_2=t - \tau_{M_A}$.
Thus, the requirement of exactness for all linear functions can be formulated
as
\begin{multline*}
w_2 \alpha_2
= \int_0^t \ell(\tau) 
   \exp[\delta_{M_A}(\tau-t)] \,\rd\tau
= \frac{1-(1+\delta_{M_A}t)\exp(-\delta_{M_A}t)}{\tau_{M_A}\delta_{M_A}^2} \alpha_2
\\
+ \frac{ (1+\delta_{M_A}(t-\tau_{M_A}))\exp(-\delta_{M_A}t) - 1 + \delta_{M_A}\tau_{M_A}}
{\tau_{M_A}\delta_{M_A}^2} \alpha_1.
\end{multline*}
This equality is satisfied for all $\alpha_1$ and $\alpha_2$ if
$$
  \tau_{M_A} = \frac{1-(1+\delta_{M_A}t)\exp(-\delta_{M_A}t)}{\delta_{M_A}(1-\exp(-\delta_{M_A}t))}
  \quad\text{and}\quad
  w_2 = \frac{1-\exp(-\delta_{M_A}t)}{\delta_{M_A}}.
$$
As above, since the exponential $\exp(-\delta_{M_A}t)$ decays rapidly towards
zero, we can simplify the expressions for $\tau_{M_A}$ and $w_2$ to	
$\tau_{M_A} = 1/\delta_{M_A}$ and $w_2 = 1/\delta_{M_A}$.

Consequently, the variable $M_A$ can be approximated as
$M_A(t) = w_2 \Phi(t-\tau_{M_A}) = (\alpha_A' + (\alpha_A - \alpha_A') D_A(t-\tau_{M_A}) ) /\delta_{M_A}$. Clearly, this is the steady state value of $M_A$ evaluated
at time delayed by $\tau_{M_A} = 1/\delta_{M_A}$.

\section{Delayed quasi steady state approximation of the circadian system}

Delayed quasi steady state assumptions derived above
are now applied to equations \ttg{eq:DA}--\ttg{eq:A}.
All these equations have the same structure, namely
\begin{equation}
  \label{eq:X}
  \rd  X(t) / \rd t = f(t) - g(t) X(t)
\end{equation}
for suitably chosen coefficients $f(t)$ and $g(t)$.
The template derivation performed above for equation \ttg{eq:MA} 
is now used for equation \ttg{eq:X} and we define
formally its time delay as $\tau_X(t) = 1/g(t)$ and its delayed
quasi steady state approximation as $X^\tau(t) = f(t - \tau_X(t))/g(t-\tau_X(t))$.

Applying this methodology to \ttg{eq:DA}--\ttg{eq:A} with \ttg{eq:subst},
we obtain the following delays and approximations:
\begin{alignat}{2}
\label{eq:dDA}
 \tau_{D_A}(t) &= \left[\theta_A + \gamma_A \widetilde A^s(R(t))\right]^{-1},
 &\quad
 D_A^\tau(t) &= D_A^s(A^\tau(t - \tau_{D_A}(t))),
\\ \label{eq:dDR}
 \tau_{D_R}(t) &= \left[\theta_R + \gamma_R \widetilde A^s(R(t))\right]^{-1},
 &\quad
 D_R^\tau(t) &= D_R^s(A^\tau(t - \tau_{D_R}(t))),
\\ \label{eq:dMA}
 \tau_{M_A} &= \delta_{M_A}^{-1},
 &\quad
 M_A^\tau(t) &= M_A^s(A^\tau(t - \tau_{M_A})),
\\  \label{eq:dMR}
 \tau_{M_R} &= \delta_{M_R}^{-1},
 &\quad
 M_R^\tau(t) &= M_R^s(A^\tau(t - \tau_{M_R})),
\\  \label{eq:dA}
 \tau_A(t) &= 
 \makebox[40mm][l]{$[\gamma_A D_A^\tau(t) + \gamma_R D_R^\tau(t) + \gamma_C R(t) + \delta_A]^{-1},$}
  &\quad
&\qquad\quad A^\tau(t) = A^s(t - \tau_A(t)),
\end{alignat}
where the definition of $A^s$ comes directly from \ttg{eq:A} with \ttg{eq:subst}
and reads as
$$
 A^s(t) = \frac{\beta_A  M_A^\tau(t) + \theta_A (1- D_A^\tau(t)) + \theta_R (1- D_R^\tau(t))}{\gamma_A  D_A^\tau(t) + \gamma_R  D_R^\tau(t) + \gamma_C  R + \delta_A}.
$$
Note that $\widetilde A^s$ was defined already in \ttg{eq:tA}.

The remaining two variables $R$ and $C$ are naturally computed by their ODEs
\ttg{eq:R}--\ttg{eq:C}, where $M_R$ and $A$ are replaced by their respective
approximations:
\begin{align}
  \label{eq:dR}
  \rd  R(t) / \rd t  &=
  \beta_R  M_R^\tau(t) - \gamma_C  A^\tau(t) R(t) + \delta_A  C(t) - \delta_R  R(t),
\\
  \label{eq:dC}
  \rd  C(t) / \rd t  &=
  \gamma_C  A^\tau(t)  R(t) - \delta_A  C(t).
\end{align}
To solve this system of delayed differential equations, we constantly extend
the initial conditions \ttg{eq:icond} to the interval $(-\infty,0]$.

System \ttg{eq:dR}--\ttg{eq:dC} with \ttg{eq:dDA}--\ttg{eq:dA} is a system
of delay differential equations with state dependent delays \cite{Kuang:1993}.
The dependence of the delays on $R$ is complicated, but it can be simplified.
Instead of variable
delays $\tau_{D_A}(t)$, $\tau_{D_R}(t)$, and $\tau_{A}(t)$ we may consider
constant delays
$\tau_{D_A}^*= 1/\theta_A$, $\tau_{D_R}^* = 1/\theta_R$, and
state dependent delay $\tau_A^*(t)=1/(\gamma_C R(t) + \delta_A)$.
The effect of this simplification is numerically assessed in the following
section.

\section{Numerical assessment of the accuracy}
\label{se:numer}

System \ttg{eq:dR}--\ttg{eq:dC} with \ttg{eq:dDA}--\ttg{eq:dA}
can be easily solved numerically.
We obtained high accuracy by implicit Euler method with time
step $10^{-3}$\,h.
Three components of this numerical solution
are presented in Figure~\ref{fi:1} (right panel) as dashed lines.

Comparing the two panels in Figure~\ref{fi:1} we clearly see that
delayed approximations are much more accurate.
However, in order to quantify the accuracy we measure and compare
both the period and amplitude of oscillations.

Let $\bff(t) = (D_A,D_A',D_R,D_R',M_A,M_R,A,R,C)(t)$ be the solution of
the original system \ttg{eq:DA}--\ttg{eq:C}. Except for an initial transient,
it is a periodic vector with period $\porig$.
Similarly, let $\bfg(t)$ be a solution of one of the approximate systems
described above and let its period be $\papprox$. To quantify the accuracy
of proposed approximations, we define the relative error in the period
and the relative $L^2$-error as 
$$
  \Erelper = \frac{|\porig - \papprox|}{\porig}
  \quad\text{and}\quad
  \ErelL = \frac{\|\bff - \widetilde\bfg\|_{L^2(a,b)}}{\|\bff\|_{L^2(a,b)}},
$$
where $b-a=\porig$ and
$\widetilde\bfg$ is linearly scaled and shifted vector $\bfg$ such that the error
in period is eliminated. In particular, $\widetilde\bfg$ is linearly scaled such that
its period is exactly $\porig$ and it is shifted such that local maxima of
$\bff$ and $\widetilde\bfg$ match.

Table~\ref{ta:1} presents periods and relative errors for the original system
\ttg{eq:DA}--\ttg{eq:C} and for its various approximations.
Namely, the third column corresponds to problem \ttg{eq:rR}--\ttg{eq:rC},
which is the original system simplified by standard quasi steady
state assumptions.
The fourth column shows problem \ttg{eq:dR}--\ttg{eq:dC} with delays \ttg{eq:dDA}--\ttg{eq:dA}
and the fifth column presents the same problem with simplified delays
$\tau_{D_A}^*$, $\tau_{D_R}^*$, and $\tau_A^*$.
The last column corresponds to the same case as the fifth, but
the only state dependent delay $\tau_A^*$ is replaced by constant
delay $\tau_A^{**} = \tau_{M_A}$. 

\begin{table}
\begin{tabular}{r|c|c|c|c|c}
 & original &  \multicolumn{4}{|c}{approximations} \\
 & system   & standard    & derived & simplified & constant \\
 &          & (no delays) & delays  & delays     & delays \\ \hline
period     & 25.6\,h & 17.9\,h & 25.1\,h & 25.3\,h & 26.1\,h \\
$\Erelper$ & --- & 29.8\,\% & 1.65\,\% & 1.02\,\% & 2.28\,\% \\
$\ErelL$   & --- & 92.7\,\% & 19.0\,\% & 19.0\,\% & 22.7\,\% \\
\end{tabular}
\bigskip
\caption{\label{ta:1}
Period and relative errors for various approximations.
}
\end{table}

We clearly observe quantitatively poor approximation properties of 
the standard quasi steady state assumptions. However, introducing delays yields
approximations with relative errors in the period about 1--2\,\% only.
The amplitude and the actual shape of the solution measured by relative
$L^2$-error considerably improved as well.


\section{Conclusions}
\label{se:conclusion}

Results presented above indicate that introducing delays to standard quasi
steady state assumptions considerably improves the quantitative accuracy
of the reduced system. Further, the derivation based on
quadrature formulas results in explicit expressions for the actual size
of the delays. The presented numerical results demonstrate the accuracy
of the proposed approach and show that the derived complicated dependence
of delays on the state variables can be simplified up to constant delays.

Since the studied model of circadian rhythms is based on standard biochemical
processes such as transcription and translation, the presented technique of
delayed quasi steady state assumptions can be easily applied
to many other biochemical networks. The simplest case is the kinetics of
mRNA, where the derived delays are constant and inversely proportional to
the decay rates $\delta_A$ and $\delta_R$.
More complicated expressions for delays were derived in the
case of genes. However, their fast kinetics and the fact that there is just
one molecule of DNA in a cell enables to simplify these delays to constants
with practically no influence on the accuracy. Finally, the complicated
dependence of delays on the state variables for some proteins
can be simplified up to constant delays.

Certain mathematical models of gene expression are based on
ODEs \cite{Vilar:2002,XieKul:2007} others are based on delay differential
equations \cite{CheAih:2002,VerRan:2008}.
The presented study may contribute to the understanding of the connection
between these models and it may suggest that models with and without delays
are just two sides of the same coin.

%

{\small
{\em Authors' addresses}:
{\em Tom\'a\v{s} Vejchodsk\'y}, University of Oxford, Mathematical Institute,
Radcliffe Observatory Quarter, Woodstock Road, Oxford, OX2 6GG, United Kingdom;
and
Institute of Mathematics, Academy of Sciences, \v{Z}itn\'a 25, Praha~1,
CZ-115\,67, Czech Republic;
e-mail: \texttt{vejchod@\allowbreak math.cas.cz}.
}

\end{document}